\documentclass[12pt,a4paper]{article}

\usepackage{amssymb}

\newtheorem{Proposition}{Proposition}[section]
\newtheorem{Theorem}[Proposition]{Theorem}
\newtheorem{Corollary}[Proposition]{Corollary}
\newtheorem{Lemma}[Proposition]{Lemma}

\newcommand {\F}{\mbox{${\cal F}$}}
\newcommand {\Fl}{\mbox{${\cal F} \!\!\!\!\!\!\!\; {\cal J}$}}
\newcommand {\Abb}[5]{\begin{eqnarray*} #1  #2  & \longrightarrow & #3 \\ 
              #4 & \longmapsto & #5 \end{eqnarray*} }
\newcommand {\aequi}{\mbox{$\Leftrightarrow$}}
\newcommand{\subn}{\mathop{\subset}\limits_{\neq}}
\newcommand{\proof}{\noindent {\bf Proof: }}
\newcommand {\qed}{\hfill $\square$}
\newcommand {\su}{\mbox{$ \subset $}}

\def\mwmatrix#1{\null\,\vcenter{\openup1\jot \mathsurround00pt\ialign{&\strut\hfil$\displaystyle{{}##{}}$\hfil\crcr#1\crcr}}\,}

\begin{document}

\title{The cohomology of period domains for reductive groups over finite fields}

\author{Sascha Orlik}
\date{}
\maketitle

\section{Introduction}
The goal of this paper is to give an explicit formula for the $\ell$-adic
cohomology of period domains over finite fields for arbitrary
reductive groups. The result is a generalisation of the computation in
\cite{O} which treats the case of the general linear group.

Let $k=\mathbb F_q$ be a finite field and $G$ a reductive algebraic group
defined over $k$ of semisimple rank  $d:=\mbox{rk}_{ss}(G)$ resp. semisimple $k$-rank $d':= k\mbox{-rk}_{ss}(G).$ 
Fix a pair $T \subset B$ consisting of a maximal torus $T$ and a Borel
subgroup $B$ both defined over $k$. We denote by $R$ the roots, by $R^+$ the positive roots
and by $$\Delta=\{\alpha_1,\ldots,\alpha_d\}$$ the base of simple roots
with respect to $T \subset B.$ 
Let  $k'$ a finite
field extension of $k$ over which $G$ splits and $\overline{k}$ an
algebraic closure of $k.$ We denote by $\Gamma=\mbox{Gal}(k'/k)$
resp. $\Gamma_k=\mbox{Gal}(\overline{k}/k)$ the associated
Galois groups. We have a natural action of $\Gamma$ on
$X^\ast(T)$ which preserves $\Delta$  since $G$ is quasi-split.
Let $$\Delta/\Gamma = \{\bar{\alpha}_1,\dots,\bar{\alpha}_{d'}\}$$ be the set of
orbits. If $\lambda \in X_\ast(T)_{\mathbb Q}$ is a rational cocharacter we will denote by 
$$P(\lambda)=\{g \in G;\, \lim_{t\to 0}\lambda(t)g\lambda(t)^{-1}
\mbox{exists in } G\}$$
%For a given root $\alpha \in R$ we denote by $U_\alpha$ the associated one-dimensional additive  root subgroup.
%&=&T\langle U_\alpha;\alpha \in R, \langle
%  \lambda,\alpha\rangle \geq 0 \rangle \\ \\ &=& T\langle U_\alpha;\alpha \in \Delta\cup -\Delta, \langle \lambda,\alpha\rangle \geq 0 \rangle$$
the associated parabolic subgroup and by 
$$U(\lambda)=\{ g\in G;\, \lim_{t \to 0}\lambda(t)g\lambda(t)^{-1} =1\}$$ its
unipotent radical.

Fix a conjugacy class $$\{\mu\} \subset X_\ast(G)$$ of one-parameter
subgroups (1-PS)  of $G$, where $\mu$ denotes a representative lying in
$X_\ast(T).$ Let 
$$E=\big\{x\in \overline{k};\, \sigma(x)=x\, \forall\, \sigma \in
Stab_{\Gamma_k}(\{\mu\})\big\}$$ be the Shimura field of $\{\mu\},$ an intermediate field of $k'/k.$ 
According to a lemma of Kottwitz (\cite{K} Lemma 1.1.3 ) there exists an
element $\mu\in\{\mu\}$ that is defined over $E.$ 
Hence the conjugacy class $\{\mu\}$ defines a flag variety
$$ \Fl:= \Fl(G,\{\mu\}):=G/P(\mu)$$ over $\overline{k}$ that is defined
over $E.$ Notice that the geometric points of $\Fl$ coincide with the set 
$$\{\mu\}/\sim \; ,$$ where $\lambda_1, \lambda_2 \in \{\mu\}$ are equivalent, written $\lambda_1 \sim 
\lambda_2,$ if there exists  $g \in U(\lambda_1)$ with $g\lambda_1g^{-1}=\lambda_2.$ Finally we set $\Gamma_E:=\mbox{Gal}(k'/E).$ 

In the further text we identify a variety with the set of its closed points.
Let $x \in \Fl$ be a point which is represented by a 1-PS
$\lambda.$ It is well known that $\lambda$ induces for every $G$-module $V$
over $\overline{k}$ a descending $\mathbb Z$-filtration
$\F^\bullet_\lambda(V)$ on $V.$ In fact, let $V=\oplus V_\lambda(i)$ be the
associated $\mathbb Z$-grading. Then $\F^\bullet_\lambda(V)$ is given by
$$ \F_\lambda^i(V)=\bigoplus_{j\geq i}V_\lambda(j), \; i \in \mathbb Z.$$
As this filtration depends only on
$x,$ we denote this filtration by $\F^\bullet_x(V).$  Considering the
adjoint action of $G$ on  its Lie algebra $Lie\;G,$ we get in particular a filtration
$\F^\bullet_x:=\F^\bullet_x(Lie\; G)$ on $Lie\; G.$ We will say that $x$ is semistable if the filtered
vectorspace $(Lie\; G, \F^\bullet_x)$ is semistable. For the latter
definition of semistability confer \cite{R1} -\cite{R3} or \cite{O} Def. 1.13. Following \cite{R3}
the semistable points of $\Fl$ form an open subvariety 
$$\Fl^{ss}:=\Fl(G,\{\mu\})^{ss},$$
which is called the period domain with respect to $G$ and $\{\mu\}.$ It is
defined over $E.$ In his paper \cite{T} Totaro has
shown that there exists a relationship to the concept of
semistability in Geometric Invariant Theory introduced by Mumford
\cite{M}. We shall explain this relationship in \S 2.

Choose an invariant
inner positive definite product on $G.$ I.e. we have for all maximal tori
$T$ in $G$ a non-degenerate positive definite pairing $(\, ,\, )$ on $X_\ast(T)_{\mathbb Q}$, such that the natural maps
$$  X_\ast(T)_{\mathbb Q} \longrightarrow X_\ast(T^g)_{\mathbb Q}$$
induced by conjugating with $g \in G(k)$ and
$$  X_\ast(T)_{\mathbb Q} \longrightarrow X_\ast(T^\sigma)_{\mathbb Q}$$
induced by conjugating with $\sigma \in \Gamma_k$
are isometries for all $g\in G(k), \sigma \in \Gamma_k.$ Here $T^g=gTg^{-1}$ is the conjugate torus resp. $T^\sigma=\sigma \cdot
T$ is the image of $T$ under the morphism $\sigma:G\rightarrow G$ induced
by $\sigma.$
The inner product, together with the natural pairing
$$ \langle \,,\, \rangle: X_\ast(T)_{\mathbb Q} \times X^\ast(T)_{\mathbb
  Q} \longrightarrow \mathbb Q \:,$$
induces identifications
\Abb{}{X_\ast(T)_{\mathbb Q}}{X^\ast(T)_{\mathbb Q}}{\lambda}{\lambda^\ast}
resp.
\Abb{}{X^\ast(T)_{\mathbb Q}}{X_\ast(T)_{\mathbb Q}}{\chi}{\chi^\ast}
 for all maximal tori $T$ in $G$. We call $\lambda^\ast$ the dual character of $\lambda$ and $\chi^\ast$ the dual cocharacter of $\chi.$
Finally we get in a similar way an invariant inner positive definite product $( \,, \,)$ on
$X^\ast(T)_{\mathbb Q}$ such that for a given root $\alpha \in R$ the
coroot $\alpha^\vee \in X_\ast(T)_{\mathbb Q}$ coincides with 
$\frac{\textstyle 2}{\textstyle (\alpha,\alpha)}\alpha^\ast.$

%\begin{Lemma} The parabolic subgroups $P(\omega^\ast_\alpha)$ induced by the dual cocharacters $\omega^\ast_\alpha$ are exactly the maximal parabolic subgroups of $G$ that contain $B.$
%\end{Lemma} 
%\proof Indeed, let $\alpha \in \Delta,\beta \in \Delta \cup -\Delta.$ Then
%$$\langle\omega^\ast_\alpha,\beta\rangle = (\omega^\ast_\alpha, \beta^\ast)
%= (\omega^\ast_\alpha, \frac{\textstyle (\beta,\beta)}{\textstyle 2}
%\beta^\vee) = \frac{\textstyle (\beta,\beta)}{\textstyle 2}(\omega_\alpha^\ast,
%\beta^\vee)= \frac{\textstyle (\beta,\beta)}{\textstyle 2}\langle
%\omega_\alpha,\beta^\vee\rangle .$$ This expression is negative exactly for
%$\beta=-\alpha.$ So $P(\omega_\alpha^*)$ is generated by $T$ and all
%$U_\beta$ s.t. $\beta \neq -\alpha.$ This is a maximal parabolic subgroup
%of $G$ that contain $B$ and every such subgroup can be written in such a way. \qed 

 Let $$\{\omega_\alpha; \alpha \in \Delta \}\subset
X_\ast(T)_{\mathbb Q}$$ be the set of fundamental weights of $G$. 
As we have already remarked at the beginning there is a natural action of
$\Gamma$ on $\Delta.$ Therefore this action induces a permutation on the set  $\{\omega_\alpha^\ast; \alpha \in \Delta\}.$ We define for every $\alpha \in \Delta$ the following $k$-rational cocharacter of $X_\ast(S)_\mathbb Q:$
$$\tilde{\omega}^\ast_\alpha:= \sum_{\sigma \in \Gamma}\sigma \omega^\ast_\alpha.$$
Notice that if two simple roots $\alpha,\beta \in \Delta$ lie in the same orbit with respect to $\Gamma$ then the cocharacters $\tilde{\omega}^\ast_\alpha$ and $\tilde{\omega}^\ast_\beta$ will coincide. 

Before we can state the main result of this paper, we have to introduce a
few more notations. Let $W_\mu$ be
the stabilizer of $\mu$ with respect to the action of $W$ on $X_\ast(T).$ 
We denote by $W^\mu$ the set of Kostant-representatives with respect to $W/W_\mu.$
Consider the action of $\Gamma_E$ on $W.$ Since $\mu$ is defined over $E,$
this action preserves $W^\mu.$ Denote the corresponding set of orbits by
$W^\mu/\Gamma_E$ and its elements by $[w],$ where $w$ is in $W^\mu.$ Clearly the length of an element in $W$ only depends on its orbit. So
the symbol $l([w])$ makes sense. For any orbit $[w]$ we set 
$$ind_{[w]}:=Ind_{Stab_{\Gamma_E}(w)}^{\:\Gamma_E}\mathbb Q_\ell. $$
This induced representation is clearly independent of the specified representative.  
For every $\alpha \in \Delta/\Gamma$ we set 
$$\tilde{\omega}_\alpha:=
(\tilde{\omega}_\alpha^\ast)^\ast \in X^\ast(T)_{\mathbb Q}.$$ If $I\subset
\Delta/\Gamma$ is any subset we define
$$\Omega_I:=\{[w] \in W^\mu/\Gamma_E;\; \langle w\mu,\tilde{\omega}_\alpha
\rangle > 0 \; \forall \alpha \not\in I\}.$$
We get the following inclusion relation
$$ I \subset J \Rightarrow \Omega_I \subset \Omega_J.$$
In the further text we denote for $[w] \in W^\mu/\Gamma_E$  by $I_{[w]}$ the 
smallest subset of $\Delta/\Gamma$ such that $[w]$ is contained in
$\Omega_{I_{[w]}}.$  Obviously we have 
\begin{equation} I_{[w]} \subset I \; \aequi \;[w] \in \Omega_I.\end{equation}

For a parabolic subgroup $P\subset G$ defined over $k$ we consider the trivial
representation of $P(k)$ on $\mathbb Q_\ell.$ We denote by 
$$i^G_P=i_{P(k)}^{G(k)}(\mathbb Q_\ell)$$ 
the resulting induced representation of $G(k).$ Further we set 
$$v^G_P=i^G_P/\sum_{P\subn Q}i^G_Q.$$ In the case $P=B$
we get the Steinberg representation \cite{S}. Finally if $I
\subset \Delta/\Gamma$ we set 
$$ P_I:=\bigcap_{I \subset \Delta/\Gamma-
  \{\alpha\}}P(\tilde{\omega}^\ast_\alpha).$$
This parabolic subgroup is defined over $k$ since the $\tilde{\omega}_\alpha^\ast$ are. Thus we can 
state the following theorem, which calculates the $\ell$-adic cohomology with
compact support of the period domain $\Fl^{ss}$ as representation of the
semi-direct product $G(k) \rtimes  \Gamma_E.$

\begin{Theorem}
We have
 $$H^{\ast}_c(\Fl^{ss},\mathbb Q_{\ell})=\bigoplus_{[w] \in W^\mu/\Gamma_E}
v^G_{P_{I_{[w]}}}\otimes
ind_{[w]}(-l([w]))[-2l([w])-\#(\Delta/\Gamma-I_{[w]})].$$ {\rm Here the symbol $(n), n\in \mathbb N$ means the $n$-th Tate twist and
$[-n], n \in \mathbb N$ symbolizes that the
corresponding module is shifted into degree $n$ of the graded cohomology ring.}
\end{Theorem}

As in the case of the $GL_{d+1}$ (cf. \cite{O} Korollar 4.5) we can state the
following result about the vanishing of some cohomology groups of period
domains. The proof of this corollary is similar to the $GL_{d+1}$-case.
\begin{Corollary}
We have
$$ H^i_c(\Fl^{ss},\mathbb Q_\ell)=0, \;\;\; 0 \leq i \leq d'-1$$
and
$$ H^{d'}_c(\Fl^{ss},\mathbb Q_\ell)= v^G_B .$$
\end{Corollary}

Theorem 1.1 has been conjectured by Kottwitz and Rapoport, who had calculated
previously the Euler-Poincar\'e characteristic with compact support of these period domains in the
Grothendieck group of $G(k) \rtimes \Gamma_E$ representations
(cf. \cite{R3}). The formula for the Euler-Poincar\'e characteristic is
accordingly $$ \chi_c(\Fl^{ss}_g,\mathbb Q_\ell)=\sum_{[w] \in
  W^\mu/\Gamma_E}(-1)^{d'-\#I_{[w]}} v^G_{P_{I_w}}\otimes ind_{[w]}(-l([w])).$$ 
In the split case the formula of the theorem  becomes
$$H^\ast_c(\Fl^{ss}, \mathbb Q_\ell)=\bigoplus_{w \in W^\mu}
v^G_{P_{I_w}}(-l(w))[-2l(w)-(\Delta -\#I_w)],$$
which has been already calculated for  $G=GL_{d+1}$ in a slightly different
way  in \cite{O}. At this point I thank Burt Totaro, Michael Rapoport and Ulrich G\"ortz for useful
discussions. Especially I am grateful for the decisive idea of Burt Totaro
in the proof of  the acylicity of the fundamental complex introduced in \S
3. His proof is based upon geometric invariant theory, whereas 
the proof in \cite{O} uses a contraction lemma of Quillen (cf. \cite{Q}
1.5). Finally I want to thank the Network for Arithmetic Algebraic Geometry of the European Community for
their financial support of my stay in Cambridge in March, where this paper was developed.

\section{The relationship of period domains to GIT}

In this section we want to explain the relationship between period
domains and Geometric Invariant Theory. For details we refer to the papers 
\cite{T} resp. \cite{R3}. We mention, that Totaro has described in his
article \cite{T} the theory of period domains in the case of local
fields. But as the reader verifies easily, all the proofs and ideas work also in the case of finite fields.

Let $$M:=P(\mu)/U(\mu)$$ be the Levi-quotient of $P(\mu)$ with center  $Z_M.$ Then $\mu$ defines an element of $X_\ast(Z_M).$ Let $T_M$ be a maximal torus in $M.$ Then we have $Z_M \subset T_M$ and $T_M$ is the isomorphic image of a maximal torus in $G.$ So we get an invariant inner product on $M$. 
Consider the dual character 
$$\mu^\ast \in X^\ast(T_M)_{\mathbb Q}.$$
As $\mu$ belongs to $X_\ast(Z_M),$ the dual character $\mu^\ast$ is contained in 
$$X^\ast(M_{ab})_{\mathbb Q} \cong Hom(P(\mu),\mathbb G _m)\otimes_{\mathbb Z} \mathbb Q.$$ The inverse character $-\mu^\ast$ induces a homogeneous line bundle $\cal{L:=L_{-\mu^\ast}}$ on $\Fl.$ The reason for the sign is that this line bundle is ample.

Let $\lambda:\mathbb G_m \rightarrow G$ be a 1-PS of $G.$ For any point $x\in
\Fl$ we can consider the slope $\mu^{\cal L}(x,\lambda)$ of $\lambda$ in $x$
relative to the line bundle ${\cal L}$ (cf. \cite{M} Def. 2.2 ). Now we are
able to state the following theorem of Totaro (cf. \cite{T} Theorem 3).

\begin{Theorem} (Totaro) 
Let $x$ be a point of $\Fl$. Then $x$ is semistable if and only if
$\mu^{\cal L}(x,\lambda) \geq 0 \;\; \forall $1-PS $\lambda$ of $G_{der}$
which are defined over $k.$ Here $G_{der}$ is the derived group of $G.$
\end{Theorem}

In order to investigate the GIT-semistability of points on varieties, it is
useful to consider the spherical building of the given group.
Let $B(G)_k$ be the real  $k$-rational spherical building of our fixed group
$G.$ Recall the definition of $B(G)_k$ (cf. \cite{CLT}). For a maximal
$k$-split torus $S$ of $G$ we consider first of all the space of rays
$$(X_\ast(S)_{\mathbb R}-\{0\}) / {\mathbb R}_{>0}:=\{\mathbb
R_{>0}\lambda; \lambda \in X_\ast(S)_{\mathbb R}-\{0\} \} $$ in
$X_\ast(S)_{\mathbb R}$ starting in the origin. This space is homeomorphic
to the $(r-1)$-sphere $\mathbb S^{r-1},$ where $r$ is the $k$-rank of $G.$
We can associate to every ray $\mathbb
R_{>0}\lambda \in$ $  (X_\ast(S)_{\mathbb R}-\{0\}) / {\mathbb R}_{> 0}$ 
a well-defined parabolic subgroup
$P(\mathbb R_{>0}\lambda)\; (\mbox{cf. }\cite{CLT}),$ 
which is compatible with the old definition of $P(\lambda)$ with respect to
a rational 1-PS $\lambda \in X_\ast(S)_{\mathbb Q}.$ 
We also have a natural action of the $k$-rational points of $G(k)$ on the
disjoint union $\coprod_{S\; k\mbox{\scriptsize -split}}(X_\ast(S)_{\mathbb
  R}-\{0\})/\mathbb R_{>0}.$ We will say that two rays $\mathbb
R_{>0}\lambda_1, \mathbb R_{>0}\lambda_2$ are equivalent, $\lambda_1\sim
\lambda_2,$ if there exists an element $g\in P(\mathbb R_{>0}\lambda_1)(k)$ which transforms the one  ray into the other. Finally we set
$$B(G)_k:= \Big(\coprod_{S\; k\mbox{\scriptsize -split}}(X_\ast(S)_{\mathbb R}-\{0\})/\mathbb
R_{>0}\Big)\Big /\sim $$ and supply this set with the induced topology. 
Again we can associate to every point $x \in B(G)_k$ a well-defined
$k$-rational parabolic subgroup $P(x)$ of $G$.
If $S$ is any maximal $k$-split torus then we have an injection
$B(S)_k \hookrightarrow B(G)_k.$ 
The space $B(S)_k$ is called the apartment belonging to $S.$ 

Assume for the remainder of this section that our group $G$ is
semisimple. In this case the space $B(G)_k$ is homeomorphic to the
geometric realisation of the combinatorial building (cf. \cite{CLT}, 6.1). Thus we have a simplicial structure on $B(G)_k$ which is defined
as follows. For a $k$-rational parabolic subgroup $P\subset G$ we let
$$D(P):=\{ x \in B(G)_k; P(x) \supset P\}$$
be the facette corresponding to $P.$ If $P$ is a minimal parabolic subgroup,
i.e. a Borel as $G$ is quasi-split, then we call $D(P)$ a chamber of
$B(G)_k.$ If in contrast $P$ is a proper maximal subgroup then $D(P)$ is called a vertex. 

Consider the $k$-rational cocharacters $\tilde{\omega}_\alpha^\ast,\alpha\in
\Delta,$ introduced in the previous section. 
%Notice that the parabolic
%subgroups $P(\tilde{\omega}^\ast_\alpha)$ are defined over $k,$ since the
% $\tilde{\omega}^\ast_\alpha$ are $k$-rational.
These cocharacters correspond then to the vertices of the chamber $D_0:=D(B),$ 
%which follows easily from the description in \cite{Ro} 2.9.
since the $P(\tilde{\omega}^\ast_\alpha), \alpha \in \Delta/\Gamma,$ are the
maximal $k$-rational parabolic subgroups that contain $B$.
For any other chamber $D=D(P)$ in $B(G)_k,$ there exists a $g\in G(k),$ such
that the conjugated elements $g\tilde{\omega}_\alpha^\ast g^{-1}, \alpha \in
\Delta/\Gamma,$ correspond to the vertices of $D.$ The element $g$ is of
course unique up to multiplication by an element of $P(k)$ from the
left. Therefore we choose for the rest of this paper for every chamber $D$ an element $g_D$ with the
above property. The element $g_{D_0}$ should be of course the obvious one. 
With this choice, we define for every chamber $D$ in $B(G)_k$ the simplex
$$\tilde{D}:=\{\sum_{\alpha \in \Delta/\Gamma} r_\alpha \lambda_\alpha\,;\:
0\leq r_\alpha\leq 1, \sum_\alpha r_\alpha=1\} \subset X_\ast(S)_{\mathbb R},$$ 
which is the convex hull of the fixed set of representatives $\lambda_\alpha:=
g_D\tilde{\omega}^\ast_\alpha g_D^{-1}  \in
X_\ast(S)_{\mathbb  R} , \alpha \in \Delta/\Gamma.$ The topological spaces $D$ and $\tilde{D}$ are obviously homeomorphic.
For the standard chamber $D_0$ we have in particular the description
$$\tilde{D}_0:=\{\sum_{\alpha \in \Delta/\Gamma} r_\alpha
\tilde{\omega}^\ast_\alpha\,;\: 0\leq r_\alpha\leq 1, \sum_{\alpha} r_\alpha=1\}.$$ 
We can extend $\mu^{\cal L}(x,\cdot) $ in a well-known way to a function on
$X_\ast(T)_{\mathbb R}$ for every maximal torus $T$ in $G.$
Notice that the slope function $\mu^{\cal L}(x,\cdot) $ is not defined on $D$ but on $\tilde{D}.$ In spite of this fact we we will say that
$\mu^{\cal L}(x,\cdot)$ is affine on $D$ if  it is affine on $\tilde{D},$ i.e.
if following equality holds:
$$ \mu^{\cal L}(x,\sum_{\alpha \in \Delta/I} r_\alpha \lambda_\alpha
)=\sum_{\alpha \in \Delta/I} r_\alpha\mu^{\cal L}(x,\lambda_\alpha) \;\;\;\; \forall \sum_{\alpha \in \Delta/I}r_a\lambda_\alpha \in \tilde{D}. $$

%\begin{Lemma}
%Let $D$ be a maximal simplex in $\Delta(G)_k.$ Then
%\begin{eqnarray*} P(\lambda_1 + \lambda_2)& = &P(\lambda_1)\cap
%  P(\lambda_2) \\
%M(\lambda_1 + \lambda_2) &=& M(\lambda_1) \cap M(\lambda_2) \; \forall
%\Lambda_1, \lambda_2 \in \tilde{D}
%\end{eqnarray*}
%\proof We may suppose  that $D=D_0$ and $G$ is split. For an arbitrary
%$\lambda \in D$ we then have
%\begin{eqnarray*}
%P(\lambda) &=& B \langle U_\alpha; \langle \lambda, \alpha \rangle \geq 0
%\; \alpha \in -\Delta \} \\
%M(\lambda) &=& T \langle U_\alpha; \langle \lambda, \alpha \rangle = 0
%\; \alpha \in \Delta \cup -\Delta \}.
%\end{eqnarray*}
%But $\langle \alpha,\lambda\rangle \geq 0 \aequi \langle
%\alpha,\lambda\rangle=0$ since $\lambda$ lies in $D_0.$
%So $\langle \alpha, \lambda_1+\lambda_2\rangle = \langle \alpha_1 \rangle +
%\langle \alpha, \lambda_2\rangle = 0 \aequi \langle \alpha,
%\lambda_1\rangle = 0 = \langle \alpha,\lambda_2\rangle.$ \qed
%\end{Lemma}

In the case of the special linear group we can calculate the slope of a
point explicitly. If $\F$ and $\F'$ are two filtrations on a finite
dimensional vector space $V$ we set
$$ (\F,\F')=\sum_{\alpha,\beta \in \mathbb Z} \alpha\beta \dim
gr^\alpha_{\F}(gr^\beta_{\F'}(V)).$$

\begin{Lemma} Let $G=SL(V).$ \vspace{-0.1cm}

\noindent (i) Let $x \in \Fl$ and $\lambda \in X_\ast(G)$ with corresponding
filtration $\F_\lambda$ on $V_{\overline{k}}=V\otimes_k\overline{k}.$ Then
$$\mu^{\cal L}(x,\lambda)=-(\F_x(V_{\overline{k}}),\F_\lambda).$$

\medskip
\noindent (ii) Let $T\subset G$ be a maximal torus and $\lambda,\lambda' \in
X_\ast(T).$
Then $$ (\lambda,\lambda')=(\F_\lambda,\F_{\lambda'}).$$
\end{Lemma}

\proof (i) If the point $x$ is fixed by $\lambda$ then our statement is
just a result of Totaro (cf. \cite{T} Lemma 6 and part (ii) of this lemma).
In general let $x_0:=\lim_{t \to 0}\lambda(t) x \in \Fl.$ Then we know that
$\mu^{\cal L}(x,\lambda)=\mu^{\cal L}(x_0,\lambda)$ (cf. \cite{M} Def. 2.2
, property (iv)). On the other hand let $\F^\bullet_x(V)$
resp. $\F^\bullet_{x_0}(V)$ be the corresponding filtrations on $V.$ Then
we claim that 
$$ gr_{\F_\lambda}^\alpha(\F_x^\beta(V)) \cong
gr_{\F_\lambda}^\alpha(\F_{x_0}^\beta(V))\: \forall \alpha, \beta \in
\mathbb Z ,$$
proving our assertion. Indeed let $W \subset V$ be any subspace. For every 
$\alpha \in \mathbb Z$ we set
$$ W_\alpha:= im \big( gr^\alpha_{\F_\lambda}(W) \hookrightarrow
V_\lambda(\alpha)\big),$$ where $V=\bigoplus_\alpha V_\lambda(\alpha)$ is
the grading of $V,$ which is induced by $\lambda.$  Then we get $$ \lim_{t \to 0}\lambda(t)\cdot W = \sum_\alpha
W_\alpha ,$$ considered as points of the corresponding grassmanian variety.
But then $$gr^\alpha_{\F_\lambda}(W) \cong W_\alpha =
gr^\alpha_{\F_\lambda}(\lim_{t \to 0} \lambda(t)\cdot W)\;\;\; \forall \alpha,$$
and the claim follows.

\noindent (ii) Choose a basis of $V$ such that $T$ is the diagonal torus of $SL_{d+1}.$
Then we may identify $\lambda$ resp. $\lambda'$ with $d+1$-tuples
$\lambda=(\lambda_1,\dots,\lambda_{d+1})$
resp. $\lambda'=(\lambda_1',\dots,\lambda_{d+1}') \in \mathbb Z^{d+1}.$
Obviously we have $gr^\alpha_{\F_\lambda}(V)=V_\lambda(\alpha)$
resp. $gr^\alpha_{\F_{\lambda'}}(V)=V_{\lambda'}(\alpha)$ and
$gr^\alpha_{\F_\lambda}(gr^\beta_{\F_{\lambda'}}(V))= V_\lambda(\alpha)\cap
V_{\lambda'}(\beta) \; \forall \alpha, \beta \in \mathbb Z.$
But then 
\begin{eqnarray*} (\F_\lambda, \F_{\lambda'}) &=& \sum_{\alpha,\beta}
  \alpha\beta \dim (V_\lambda(\alpha) \cap V_{\lambda'}(\beta))\\ &=&
\sum_{\alpha,\beta}\alpha\beta \# \{ i;\; \lambda_i =\alpha, \lambda_i'
=\beta\} = \sum_{i=1}^{d+1} \lambda_i\lambda'_i =(\lambda,
\lambda'). \;\;\;\;\;\;\;\;\;\;\; \mbox{ \qed}
\end{eqnarray*}

The idea of the next proposition is due to Burt Totaro which is a decisive
point in proving the acyclicity of the fundamental complex in Theorem 3.1.

\begin{Proposition} Let $x \in \Fl$ be any point. The slope function $\mu^{\cal
    L}(x,\,\cdot )$ is affine on each chamber of $B(G)_k.$
\end{Proposition}

\proof %Since period domains only depend of the semisimple part of a given
%  reductive group (cf. Theorem 2.1), we can assume that $G=G_{der}.$
%  Further 
We may assume that our group is $k$-simple. 
 Choose a faithful representation  $$i:G \hookrightarrow SL_n=:G'$$ which is defined over $k$. Set $\mu':=i\circ \mu \in X_\ast(G').$ We get a closed immersion
$$i: \Fl(G,\{\mu\}) \hookrightarrow \Fl(G',\{\mu'\})=:\Fl'$$
of the corresponding flag varieties,
under which $\mu$ is mapped to $\mu'.$ We assume that
we have an invariant inner product on $SL_n$ which restricts to our fixed
one on $G.$ This is not really a restriction since two such inner products
on a $k$-simple group
differ only by a positive scalar (cf. \cite{T} Lemma 7). The line bundle $\cal{L':=L_{-\mu'}}$ on $\Fl',$ defined in a
similar way as $\cal{L},$  restricts then via  the pullback to $\cal{L}.$ 
 Because of the equality $\mu^{\cal L}(x,\lambda)=\mu^{\cal
  L'}(i(x),i\circ \lambda)$ (\cite{M} property (iii) following Def. 2.2) we can
  restrict ourselves to the case $G=SL(V).$ Let $\lambda \in \{\mu\}$ be a
  1-PS representing $x.$ Let $S\subset G$ be a maximal $k$-split torus, such
  that the corresponding apartement contains both $D,$ the chamber
  with representatives $\lambda_\alpha, \alpha \in \Delta/\Gamma,$ of its
  vertices  and $\lambda.$
  Using the previous lemma we get
$$ \mu^{\cal L}(x,\sum_{\alpha \in \Delta/I} r_\alpha \lambda_\alpha
)=-(\F_x(V_{\overline{k}}),\F_{\sum\limits_{\alpha \in \Delta/I} r_\alpha
  \lambda_\alpha} ) = -
(\lambda,\sum_{\alpha\in \Delta/I} r_\alpha \lambda_\alpha) = $$
$$ -\sum_{\alpha \in \Delta/I} r_\alpha (\lambda,\lambda_\alpha) =
-\sum_{\alpha \in \Delta/I} r_\alpha(\F_x(V),\F_\lambda) = \sum_{\alpha
  \in \Delta/I} r_\alpha\mu^{\cal L}(x,\lambda_\alpha).\;\;\; \square$$

%We
%  have further the equality $\mu^{\cal L}(x,\lambda) = \mu^{\cal
%  L}(gx,g\lambda g^{-1})\; \forall g \in G(k)$ \cite{M} property (i)
%  following def.2.2.  Therefore we
%  can assume that both $\lambda_x$ and all the 1-PS in question are lying
%  in the apartement defined by $T.$ But in this case we know that the slope
%  of $x$ with respect to any 1-PS $\lambda$ is nothing else than $-(\lambda_x, \lambda)$ \cite{T} lemma 6. But the latter expression is linear, so the proof is complete.

I want to stress that the previous corollary fails for arbitrary
varieties. In general the slope function is only convex (cf. \cite{M} Cor. 2.13).

\begin{Corollary} Let $x$ be a point in $\Fl.$ Then $x$ is not semistable
  $\Leftrightarrow$ There exists a $g \in G(k)$ and an  $\alpha \in \Delta$ such that $\mu^{\cal L}(x,g\tilde{\omega}^\ast_\alpha g^{-1}) < 0 .$
\end{Corollary}

\proof The direction $``\Leftarrow ``$ is clear. So let $\lambda$ be a
$k$-rational 1-PS with $\mu^{\cal L}(x,\lambda)<0.$ Let $g \in G(k)$ such
that $Int(g^{-1})\circ \lambda$ lies in the simplex $\tilde{D}_0$ spanned
by the rational 1-PS  $\tilde{\omega}^\ast_\alpha, \alpha \in \Delta/\Gamma.$ Thus we can write
$\lambda$ in the form $\lambda= \sum_\alpha r_\alpha g\tilde{\omega}^\ast_\alpha g^{-1},$ with $0 \leq r_\alpha \leq 1.$ The statement follows now immediately from Proposition 2.3.
\qed

\section{The fundamental complex}
Let $G$ be again an arbitrary reductive group.
In this section we will construct an acyclic complex of \'etale sheaves
on the closed complement 
$$Y:=\Fl \setminus \Fl^{ss}$$ of the period domain $\Fl^{ss}$ which is defined
over $E$ as well. This complex yields a
method to calculate the cohomology of $\Fl^{ss}.$

For any subset  $I \subset \Delta/\Gamma$ we set
$$ Y_I:=\{x \in \Fl;  \mu^{\cal L}(x,\tilde{\omega}^\ast_\alpha) < 0
\; \forall \alpha \notin I \}.$$

\begin{Lemma} a) The set $Y_I$ induces a closed subvariety of $Y$ defined
  over $E.$ 

\noindent b) The natural action of $G$ on $\Fl$ restricts to an action of $P_I$ on
  $Y_I$ for every $I$.
\end{Lemma}

\proof It is enough to show the statement in the extreme case
$I=\Delta/\Gamma-\{\alpha\}$ Choosing an $G$-linearized embedding $Y \hookrightarrow \mathbb
P(V)$ into some projective space (cf. \cite{M} Prop. 1.7), we may restrict
ourselves to show that the set
$ \{x\in \mathbb P(V); \mu(x,\lambda)\leq 0\} $ is closed for every $\lambda \in
X(G)_k.$ Let $V=\oplus_{i\in \mathbb Z} V(i)$ be the grading induced by
$\lambda.$ Then the above set is just the closed subspace $\mathbb
P(\oplus_{i\geq 0}V(i)),$ and the first assertion follows. The second
statement results immediately from the fact that
$$ \mu^{\cal L}(px,\tilde{\omega}_\alpha^\ast)= \mu^{\cal
  L}(x,\tilde{\omega}_\alpha^\ast)\; \forall p\in
P(\tilde{\omega}_\alpha^\ast) \; (\mbox{cf. }\cite{M}\; \mbox{Prop. 2.7} ) $$
\qed

It is a consequence of Corollary 2.4 that we can write $Y$ as the  union
$$ Y=\bigcup_{\alpha \in \Delta/\Gamma}\bigcup_{g \in G(k)} gY_{\Delta/\Gamma- \{\alpha\}}.$$
Now, let $F$ be an \'etale sheaf on $Y.$ Let $I \subset J$ be two subsets of
$\Delta/\Gamma$ with $\#(J\setminus I)=1.$ Let further $g \in (G/P_I)(k),
h\in (G/P_J)(k)$ two elements, such that $g$ is mapped to $h$  under the canonical projection
$(G/P_I)(k) \longrightarrow (G/P_J)(k).$ In this case we let
$$p_{I,J}^{g,h}:  (h\phi_J)_*(h\phi_J)^* F \longrightarrow 
(g\phi_I)_*(g\phi_I)^* F$$ be the natural morphism of \'etale sheaves on $Y$ which is induced by the
closed embedding $gY_I \hookrightarrow hY_J.$
If  $g$ is not mapped to $h$ then we set 
$p_{I,J}^{g,h}=0.$ Finally we define $$ p_{I,J}= \!\!\!\!\!\!\!\!
\bigoplus_{(g,h) \in (G/P_I)(k)\times (G/P_J)(k)}\!\!\!\!\! p_{I,J}^{g,h}: \bigoplus_{h\in (G/P_J)(k)}\!\!\!\!\!(h\phi_J)_*(h\phi_J)^* F \longrightarrow \!\!\!\!\!
\bigoplus_{g \in (G/P_I)(k)}\!\!\!\!\! (g\phi_I)_*(g\phi_I)^* F.$$ For two
arbitrary subsets  $I,J \subset S$ with $\#J - \#I=1$ we set

\[ d_{I,J} =\left\{ \begin{array}{l@{\quad : \quad }r}  (-1)^i p_{I,J} &
 J= I \cup \{\bar{\alpha}_i\} \\  0 & I \not\subset J \end{array} \right. .  \]
We get a complex of \'etale sheaves on $Y:$

$$(*): 0 \rightarrow F \rightarrow\!\!\! \bigoplus_{I \subset \Delta/\Gamma \atop
\#(\Delta/\Gamma - I)=1} \bigoplus_{g \in (G/P_I)(k)} (\phi_{g,I})_*(\phi_{g,I})^* F \rightarrow
\!\!\!\bigoplus_{I \subset \Delta/\Gamma \atop \#(\Delta/\Gamma - I)=2} \bigoplus_{g \in (G/P_I)(k)}
(\phi_{g,I})_*(\phi_{g,I})^* F \rightarrow $$
$$ \dots \rightarrow \!\!\!\bigoplus_{I \subset \Delta/\Gamma \atop \#(\Delta/\Gamma - I)=d'-1} \bigoplus_{g \in (G/P_I)(k)} (\phi_{g,I})_*(\phi_{g,I})^* F \rightarrow
\bigoplus_{g \in (G/B)(k)} (\phi_{g,\emptyset})_*(\phi_{g,\emptyset})^* F \rightarrow 0,$$ where $\phi_{g,I}$ denotes the closed immersion $gY_I \hookrightarrow Y.$

One essential step in order to calculate the cohomology of our period
domain is the following result.
\begin{Theorem}
The above complex is acyclic.
\end{Theorem}

\proof Let $x\in Y(k^{sep})$ be a geometric point. Localizing in $x$ yields a chain complex which is precisely the chain complex that computes the homology with coefficient group $F_x$ of a subcomplex of the combinatorial Tits complex to $G(k).$ Strictly speaking this subcomplex corresponds to the following subset of the set of vertices of the Tits building:
$$\{ gP(\tilde{\omega}_\alpha^\ast)g^{-1}; g \in G(k), \alpha \in \Delta
\mbox{  s.t.  } \mu^{\cal L}(x,g\tilde{\omega}^\ast_\alpha g^{-1}) <0 \}.$$
We will show that this combinatorial subcomplex is contractible. Let $T_x$
be its canonical geometric realisation in the real spherical building
$B(G)_k.$ Then $T_x$ is already contained in $B(G_{der})_k\subset B(G)_k.$ The next two lemmas will show that the topological space $T_x$ is contractible.\qed

\begin{Lemma} Let $C_x:=\{\lambda \in B(G_{der})_k ;
  \frac{\textstyle\mu^{\cal L}(x,\lambda)}{\textstyle\|\lambda\|} <0\}. $
  This set is convex. The intersection of $C_x$ with each chamber in
  $B(G_{der})_k$ is convex. ( For the definition of convex we refer to \cite{M}.) 
\end{Lemma}

\proof In the case that the group $G_{der}$ is split this is just  \cite{M}
Cor. 2.16. But the proof for the general case goes through in the same way.  \qed

\bigskip
Notice that  we get an inclusion $T_x \hookrightarrow C_x$ because the
slope-function is affine on every chamber of $B(G_{der})_k.$

\begin{Lemma} The inclusion $T_x \hookrightarrow C_x$ is a deformation retract. \end{Lemma}

\proof Let $D=gD_0, g\in G(k)$ be a chamber in the real spherical building
$B(G)_k$ with $D \cap C_x \neq \emptyset.$ Following Lemma 3.3
this intersection is a convex set, where $D \cap T_x$ lies in the boundary
of this space. Thus we can construct a deformation retract
between $D\cap C_x$ and $D \cap T_x$ in the following way.  Let $\Lambda:=
\{ \alpha \in \Delta/\Gamma : g\tilde{\omega}_\alpha^\ast g^{-1} \in T_x\}.$
  Denote by $\widetilde{D\cap C_x}$ resp. $\widetilde{D\cap T_x}$ the preimage of  $D\cap C_x$ resp. $D\cap T_x$ under the canonical homoeomorphism
  $\tilde{D} \rightarrow D.$ Let 
$$\phi_D: (D\cap C_x) \times [0,1] \longrightarrow T_x$$ be the map which is induced
by the map 
$$\phi_{\tilde{D}}: (\widetilde{D\cap C_x})\times [0,1] \longrightarrow \widetilde{D\cap T_x}$$ defined by
$$ \phi_{\tilde{D}}(\sum_{\alpha \in \Lambda} r_\alpha g\tilde{\omega}_\alpha^\ast g^{-1} + \sum_{\alpha \not\in \Lambda} r_\alpha g\tilde{\omega}_\alpha^\ast g^{-1},t):= \sum_{\alpha \in \Lambda} r_\alpha g\tilde{\omega}_\alpha^\ast g^{-1}+ \sum_{\alpha \notin \Lambda} tr_\alpha g\tilde{\omega}_\alpha^\ast g^{-1}.$$
This is a continuous map and one checks easily that the collection of these maps paste together to a continuous map
$$ \phi : C_x \times [0,1] \longrightarrow T_x$$
which induces a deformation retraction from $T_x$ to $C_x$.
\qed

\section{The proof of Theorem 1.1}

This last part of the paper deals with the evaluation of the complex $(*)$
in the case of the $\ell$-adic sheaf $F=\mathbb Q_{\ell}.$

\begin{Proposition} We have  the following description of the closed
  varieties  $Y_I$ in terms of the Bruhat cells of $G$ with respect to $P(\mu).$
\begin{eqnarray*} Y_I &=& \bigcup_{ w \in \Omega_I}
  P(\tilde{\omega}_\alpha^\ast)w P(\mu)/P(\mu) \\ &=& \bigcup_{ w \in \Omega_I} BwP(\mu)/P(\mu).
\end{eqnarray*}
\end{Proposition}

\proof It is enough to show the assertion in the case $I=\Delta/\Gamma -
\{\alpha\}$ for an element $\alpha \in \Delta/\Gamma,$ since the sets $\Omega_I$ and $Y_I$  are compatible with
forming intersections relative to the sets $I,$ i.e. 
$$ \Omega_{I\cap J} = \Omega_I \cap \Omega_J$$ 
resp.
$$ Y_{I\cap J}= Y_I \cap Y_J,\;  \forall I,J \subset \Delta/\Gamma.$$  
Let $p$ be an element of $P(\tilde{\omega}_\alpha^\ast).$ We have  the equality
$$ \mu^{\cal L}(px,\tilde{\omega}_\alpha^\ast)=\mu^{\cal
  L}(x,\tilde{\omega}_\alpha^\ast) \; \forall x \in \Fl$$ (cf. \cite{M} Prop. 2.7).
The proposition follows  now immediately from the equalities
$$  \mu^{\cal L}(pw\mu,\tilde{\omega}^\ast_\alpha)=  \mu^{\cal
  L}(w\mu,\tilde{\omega}^\ast_\alpha)=-(w\mu,\tilde{\omega}^\ast_\alpha)=
-\langle w\mu, \tilde{\omega}_\alpha \rangle. \;\;\;\;\;\; \square$$

The above cell decomposition for the varieties $Y_I$ allows us to calculate
the cohomology of them. The proof is the same as in the case of $GL_{d+1}$
(cf. \cite{O} Prop. 7.1) and will be omitted.

\begin{Proposition} \begin{equation} H_{\acute{e}t}^\ast(Y_I,\mathbb Q_\ell)=\bigoplus_{[w] \in \Omega_I}  ind_{[w]}(-l[w])[-2l([w])]\end{equation}
\end{Proposition}

%\proof The proof is the same as in the case of $GL_{d+1}$ (\cite{O}
%Prop. 7.1). In order to apply Lemma 7.2 in \cite{O},  we have merely to
%show that the corresponding index set is closed with respect to the
%Bruhat order i.e. if $w'<w$ and  $[w] \in \Omega_I$ then $[w']\in
%\Omega_I.$ But $w'<w$ implies that $w^{-1}\tilde{\omega}_\alpha$ is the sum of
%$w'^{-1}\tilde{\omega}_\alpha$ and some positive roots \cite{?} i.e.
%$$w^{-1}\tilde{\omega}_\alpha = w'^{-1}\tilde{\omega}_\alpha + \sum_{\beta
%  \in R^+}  n_\beta\beta,$$
%where $\beta$ lies in $R^+$ and $n_\beta \geq 0 \;\forall \beta.$ As $\mu$
%lies in the closure of the positive Weyl chamber the assertion follows.
%\qed

In the following we denote for an orbit $[w] \in W^\mu/\Gamma_E$ and a
subset  $I \subset \Delta/\Gamma$ the contribution of
$[w]$ with respect to the direct sum (2)  by $H(Y_I,[w])$ i.e.
 
\begin{equation}
H(Y_I,[w])=\left\{ 
\begin{array}{r@{\quad:\quad}l} ind_{[w]}(-l([w]))[-2l([w])] & [w] \in \Omega_I \\ 0 & [w]\notin \Omega_I \end{array}
\right. .
\end{equation}
Thus we have
\begin{equation}
 H^\ast_{\acute{e}t}(Y_I,\mathbb Q_\ell)= \bigoplus_{[w]\in W^\mu/\Gamma_E} H(Y_I,[w]).
\end{equation}
Let $I\subset J $ be two subsets of $\Delta/\Gamma.$ 
We consider the homomorphism
$$\phi_{I,J}: H^{\ast}_{\acute{e}t}(Y_J) \longrightarrow
H^{\ast}_{\acute{e}t}(Y_I)$$ 
given by the closed embedding $Y_I \hookrightarrow Y_J.$
The construction of Prop. 4.2 induces a grading of $\phi_{I,J},$
$$\phi_{I,J}=\!\!\!\!\!\!\bigoplus_{([w],[w']) \in W^\mu/\Gamma_E \times W^\mu/\Gamma_E}\!\!\!\!\! \phi_{[w],[w']}:
\bigoplus_{[w]\in W^\mu/\Gamma_E} H(Y_J,[w]) \longrightarrow
\bigoplus_{[w']\in W^\mu/\Gamma_E} H(Y_I,[w'])$$
with
\begin{equation} \phi_{[w],[w']}=\left\{\begin{array}{r@{\quad:\quad}l} id  & [w]=[w']
    \\ 0 & [w] \neq [w'] \end{array} \right. . 
\end{equation}

\bigskip
We need a generalisation of a result of  Lehrer resp. Bj\"orner. 
We will construct a complex in analogy to the sequence $(*).$
Let  $I\subset J
\subset \Delta/\Gamma$ two subsets with $\#(J\setminus I) =1.$ We get 
a homomorphism 
$$p_{I,J}:i^G_{P_J} \longrightarrow i^G_{P_I},$$ which comes from the
projection  $(G/P_I)(k) \longrightarrow (G/P_J)(k).$
For two arbitrary subsets  $I,J \subset \Delta/\Gamma$ with $\#J -\#I=1$
we define 
$$d_{I,J}=\left\{ \matrix{ (-1)^i p_{I,J} & J = I  \cup \{\bar{\alpha}_i\} \cr 0 & I \not\subset J } \right. .$$ 
Thus we get for every $I_0 \subset \Delta/\Gamma$ a $\mathbb
Z$-indexed complex
$$K_{I_0}^\bullet: 0 \rightarrow i^G_G \rightarrow
\!\!\!\!  \bigoplus_{I_0 \subset I \subset \Delta/\Gamma \atop
  \#(\Delta/\Gamma- I)=1}\!\!\!\! i^G_{P_I}  \rightarrow \!\!\!\!
\bigoplus_{I_0 \subset I \subset \Delta/\Gamma \atop \#(\Delta/\Gamma- I)=2}\!\!\!\! i^G_{P_I}\rightarrow
\dots \rightarrow \!\!\!\! \bigoplus_{I_0 \subset I \subset \Delta/\Gamma \atop
  \#(\Delta/\Gamma- I)=\#(\Delta/\Gamma- I_0)-1 }\!\!\!\!\!\!\! i^G_{P_I}  \rightarrow i^G_{P_{I_0}},$$
where the differentials are induced by the above $d_{I,J}.$ The component
$i^G_G$ is in degree $-1.$ 

\begin{Proposition}  The complex $K_{I_0}^\bullet$ is acyclic.
\end{Proposition}   

\proof In the split case this is precisely the result of Lehrer \cite{L} resp. Bj\"orner
\cite{Bj}. Since the group $\Gamma$ is finite taking the
fix-vectors in the category of $\mathbb Q_\ell$-representations yields an
exact functor. But the above complex is just the resulting fix-complex of the
analogous complex relative to $G$ considered as a split group defined over $k'.$ \qed

We mention the following well known lemma (cf. \cite{O}, Lemma 7.4).
\begin{Lemma} Every extension of the $Gal(\overline{k}/E)$-module $\mathbb
  Q_\ell(m)$ by $\mathbb Q_\ell(n)$ with $m\neq n$ splits.
\end{Lemma}

%\proof \cite{O} Lemma 7.4 \qed

The acyclic complex $(*)$ yields a method to calculate the $\ell$-adic
cohomology of $Y.$
% Set $$\overline{Y}:=Y\times _k \overline{k} \mbox{ bzw. }
% \overline{\Phi_{g,I}}=\Phi_{g,I} \times _k id_{\overline{k}}.$$

\begin{Theorem}   
The spectral sequence   
$$ E_1^{p,q} = H^q_{\acute{e}t}(Y,\bigoplus_{I \su \Delta/\Gamma \atop
\#(\Delta/\Gamma- I)=p+1} \bigoplus_{g \in
(G/P_I)(k)}(\phi_{g,I})_*(\phi_{g,I})^* \mathbb Q_\ell) \Longrightarrow
H^{p+q}_{\acute{e}t}(Y,\mathbb Q_\ell)$$
resulting from $(*)$
, degenerates in the $E_2$-term
and we get for the $\ell$-adic cohomology of $Y:$  
\begin{eqnarray*}H^{\ast}_{\acute{e}t}(Y,\mathbb
  Q_\ell)&=&\bigoplus_{w\in W^\mu/\Gamma_E \atop \#(\Delta/\Gamma- I_{[w]})=1} \Big(i_{P_{I_{[w]}}}^G\otimes ind_{[w]}(-l(w))[-2l([w])]\Big) \oplus \\
& &\bigoplus_{w\in W^\mu/\Gamma_E \atop \#(\Delta/\Gamma- I_{[w]})>1} 
\bigg(\Big(i^G_G \otimes ind_{[w]}(-l([w]))[-2l([w])]\Big) \oplus  \\
& & \oplus \Big(v^G_{P_{I_{[w]}}} \otimes
ind_{[w]}(-l([w]))[-2l([w])-\#(\Delta/\Gamma- I_{[w]})+1]\Big)\bigg).
\end{eqnarray*}
\end{Theorem}

\proof We have
\begin{eqnarray*} 
E_1^{p,q} &=& H^q_{\acute{e}t}(Y,\!\!\!\bigoplus_{I \subset \Delta/\Gamma \atop
\#(\Delta/\Gamma- I)=p+1 } \bigoplus_{g \in (G/P_I)(k)}(\phi_{g,I})_*(\phi_{g,I})^* \mathbb Q_\ell) \\
&=& \bigoplus_{I \subset \Delta/\Gamma \atop \#(\Delta/\Gamma- I)=p+1 }
\bigoplus_{g \in (G/P_I)(k)}H^q_{\acute{e}t}(Y_I,(\phi_{g,I})^* \mathbb
   Q_\ell) =\bigoplus_{I \subset \Delta/\Gamma \atop \#(\Delta/\Gamma-
     I)=p+1} \bigoplus_{(G/P_I)(k)} H_{\acute{e}t}^q(Y_I, \mathbb Q_\ell).
\end{eqnarray*}
The application of $(4)$ and $(5)$ yields a decomposition 
$$ E_1=\bigoplus_{[w] \in W^\mu/\Gamma_E} E_{1,[w]}$$ 
into subcomplexes with 
$$E_{1,[w]}^{p,q}=\left\{ \mwmatrix{\bigoplus_{I \subset \Delta/\Gamma \atop
      \#(\Delta/\Gamma- I)=p+1} \bigoplus_{(G/P_I)(k)} H(Y_I,[w]) & \;\;\;\;\; q=2l([w]) \cr
    0 & \;\;\;\;\; q\neq 2l([w]) } \right. .$$

\noindent Thus  $E_{1,[w]}$ is the subcomplex

$$ E_{1,[w]}: \bigoplus_{I \subset \Delta/\Gamma \atop \#(\Delta/\Gamma- I)=1} \bigoplus_{(G/P_I)(k)}H(Y_I,[w])\longrightarrow \bigoplus_{I \subset \Delta/\Gamma
\atop \#(\Delta/\Gamma- I)=2} \bigoplus_{(G/P_I)(k)}H(Y_I,[w]) \longrightarrow $$
$$ \dots \longrightarrow\bigoplus_{(G/B)(k)}H(Y_\emptyset,[w]).$$
\noindent In view of $(1)$ und $(3)$ we have
\[ H(Y_I,[w])=\left\{ 
\begin{array}{r@{\quad:\quad}l} ind_{[w]}(-l([w]))[-2l([w])] &  I_{[w]} \su I
\\ 0 & I_{[w]} \not\subset I \end{array}
\right. . \]

\noindent So $E_{1,[w]}$ simplifies to  
$$\Big(\bigoplus_{I_{[w]} \subset I \atop \#(\Delta/\Gamma-
    I)=1}i^G_{P_I}\otimes ind_{[w]}(-l([w]))
\rightarrow \bigoplus_{I_{[w]} \subset I \atop \#(\Delta/\Gamma- I) =2}
i^G_{P_I}\otimes ind_{[w]}(-l([w])) \rightarrow \dots  $$ 
$$\rightarrow i^G_{P_{I_{[w]}}}\otimes ind_{[w]}(-l([w]))\Big)[-2l([w])],$$ and we get an exact sequence of complexes:

$$ 0 \rightarrow i^G_G\otimes ind_{[w]}(-l([w]))[-2l([w])+1] \rightarrow
K_{I_{[w]}}^\bullet\otimes ind_{[w]}
(-l([w]))[-2l([w])]\rightarrow E_{1,[w]} \rightarrow 0.$$

\noindent This yields the following three cases for $E_{2,[w]}$:

\medskip
$\begin{array}{rclrl} I_{[w]}=\Delta/\Gamma & : & E_{2,[w]}^{p,q} & =& 0 \;\;  p\geq 0, q\geq 0 \\ \\
\#(\Delta/\Gamma- I_{[w]})=1 & : & E_{2,[w]}^{0,2l([w])} &=&
i^G_{P_{I_{[w]}}}\otimes ind_{[w]}(-l([w])) \\ \\
           &   & E_{2,[w]}^{p,q} &=&  0 \;\; (p,q)\neq(0,2l([w])) \\ \\
\#(\Delta/\Gamma- I_{[w]}) > 1 & : & E_{2,[w]}^{0,2l([w])} & = &
i^G_G\otimes ind_{[w]}(-l([w])) \\ \\
& & E_{2,[w]}^{j,2l([w])} &=&  0 \;\; j=1,\ldots,\#(\Delta/\Gamma-I_{[w]})-2 \\ \\
& & E_{2,[w]}^{j,2l([w]))} &=& v^G_{P_{I_{[w]}}}\otimes ind_{[w]}(-l([w])) \;\; j=\#(\Delta/\Gamma-I_{[w]})-1\\ \\
& & E_{2,[w]}^{p,q} &=& 0 \;\;  q\neq 2l([w]) \mbox{ or }  p>\#(\Delta/\Gamma- I_{[w]})-1.
\end{array}$

\medskip
\noindent The Galois modules  $E_2^{p,q}\neq (0)$ posess the Tate twist
$-\frac{q}{2}$. As every homomorphism of
Galois modules of different Tate twists is trivial, the
$E_2$-term coincides with the
$E_\infty$-term. Thus for all $n\geq 0$

$$gr^p(H_{\acute{e}t}^n(Y))=E_\infty^{p,n-p}=E_2^{p,n-p}=\bigoplus_{[w]\in W^\mu/\Gamma_E} E_{2,[w]}^{p,n-p}$$

\[ =\left\{ 
\begin{array}{r@{\quad:\quad}l}
  \bigoplus\limits_{[w] \in W^\mu/\Gamma_E \atop { \#(\Delta/\Gamma-
    I_{[w]})=1 \atop 2l([w])=n }}i^G_{P_{I_{[w]}}}\otimes
    ind_{[w]}(-l([w]))\oplus \bigoplus\limits_{[w] \in W^\mu/\Gamma_E \atop
    { \#(\Delta/\Gamma - I_{[w]})>1 \atop 2l([w])=n}} i^G_G\otimes ind_{[w]}(-l([w])) &  p=0 \\
\bigoplus\limits_{[w] \in W^\mu/\Gamma_E \atop 2l([w]) + \#(\Delta/\Gamma-
    I_{[w]})-1=n}v^G_{P_{I_{[w]}}}\otimes ind_{[w]}(-l([w]))  & p>0 \\
\end{array}
\right.  \] 

\noindent Following Lemma 4.4  extensions of $\mathbb
Q_\ell(m)$ by $\mathbb Q_\ell(n)$ with $m\neq n$ are trivial.
This yields an isomorphism
$$ H_{\acute{e}t}^n(Y,\mathbb Q_\ell )\cong \bigoplus_{p \in
  \mathbb N} gr^p(H_{\acute{e}t}^n(Y,\mathbb Q_\ell))$$
\[=\bigoplus_{[w]\in W^\mu/\Gamma_E \atop { \#(\Delta/\Gamma- I_{[w]})=1 \atop 2l([w])=n}}\!\!\!\!\!i^G_{P_{I_{[w]}}}\otimes
\:ind_{[w]}(-l([w]))\oplus\!\!\!\!\!\! \\  \bigoplus_{[w]\in W^\mu/\Gamma_E
  \atop {\#(\Delta/\Gamma-I_{[w]}) >1 \atop
  2l([w])=n}}\!\!\!\!\!i^G_G\otimes ind_{[w]}(-l([w])) \oplus \]
\[\bigoplus_{[w] \in W^\mu/\Gamma_E \atop 2l([w]) + \#(\Delta/\Gamma-
  I_{[w]})-1=n}\!\!\!\!\! v^G_{P_{I_{[w]}}}\otimes ind_{[w]}(-l([w])) .\]

\noindent The claim follows.
\qed

\bigskip

\noindent {\bf Proof of Theorem 1.1:} The proof is the same as in the case of
$G=GL_{d+1}$ (cf. \cite{O}). \qed

\end{document}